\documentclass[11pt,a4paper]{article}

\usepackage[leqno]{amsmath}
\usepackage{amsthm}
\usepackage{enumerate,graphicx}

\newcommand{\bt}{\begin{tabular}}
\newcommand{\et}{\end{tabular}}
\newcommand{\beq}[1]{\begin{equation}\label{#1}}
\newcommand{\eeq}{\end{equation}}
\newcommand\toe{\stackrel{e}{\to}}

\DeclareMathOperator\cli{cl}
\newcommand{\oneaddr}[1]{\par\small\bigskip\slshape%
                               {\tabcolsep 0pt\hfill\bt{c}#1\et}\par}
\newcommand\email[1]{e-mail\/{\normalfont: \texttt{#1}}}
\makeatletter
\renewcommand\section{\@startsection{section}{1}{\z@}%
                                    {0pt}
                                    {-0.5em}%
                                    {\normalfont\normalsize\bfseries\boldmath}}
\def\@seccntformat#1{\indent\csname the#1\endcsname.~}
\renewenvironment{proof}[1][\proofname]{\par
  \pushQED{\qed}%
  \normalfont 
        {\bfseries
    #1\@addpunct{.} }
    \ignorespaces
}{%
  \popQED
  \@endpefalse
}
\makeatother
\numberwithin{equation}{section}

\newtheorem{thm}{Theorem}[section]
\newtheorem{lemma}{Lemma}[section]
\newtheorem{pro}{Proposition}[section]
\theoremstyle{definition}
\newtheorem{definition}{Definition}[section]
\newtheorem{rem}{Remark}[section]
\begin{document}
\title{CHROMATIC NUMBER OF GRAPHS\\AND EDGE FOLKMAN NUMBERS%
\thanks{This work was supported by the Scientific Research Fund of
the St. Kliment Ohridski Sofia University under contract No~75, 2009.}}
\author{Nedyalko Dimov Nenov}

\maketitle

\begin{abstract}
We consider only simple graphs. The graph $G_1+G_2$ consists of vertex
disjoint copies of $G_1$ and $G_2$ and all possible edges between the
vertices of $G_1$ and $G_2$. The chromatic number of the graph $G$ will
be denoted by $\chi(G)$ and the clique number of $G$ by $\cli(G)$. The
graphs $G$ for which $\chi(G)-\cli(G)\ge 3$ are considered. For these
graphs the inequality $|V(G)|\ge\chi(G)+6$ was proved in~\cite{12},
where $V(G)$ is the vertex set of $G$. In this paper we prove that
equality $|V(G)|=\chi(G)+6$ can be achieved only for the graphs
$K_{\chi(G)-7}+Q$, $\chi(G)\ge 7$ and $K_{\chi(G)-9}+C_5+C_5+C_5$,
$\chi(G)\ge 9$, where graph $Q$ is given on Fig.~1 and $K_n$ and
$C_5$ are complete graph on $n$ vertices and simple 5-cycle,
respectively. With the help of this result we prove some new facts for
some edge Folkman numbers (Theorem~4.2). 


\end{abstract}

\section{Introduction.}
We consider only finite, non-oriented graphs without loops and multiple edges.
We call a $p$-clique of the graph $G$ a set of $p$ vertices each two of
which are adjacent. The largest positive integer $p$ such that $G$ contains
a $p$-clique is denoted by $\cli(G)$ (clique number of $G$).
We shall use also the following notations:
\begin{itemize}
\item
$V(G)$ is the vertex set of $G$;
\item
$E(G)$ is the edge set of $G$;
\item
$\overline G$ is the complement of $G$;
\item
$G-V$, $V\subseteq V(G)$ is the subgraph of $G$ induced by
$V(G)\setminus V$;
\item
$\alpha(G)$ is the vertex independence number of $G$;
\item
$\chi(G)$ is the chromatic number of $G$;
\item
$f(G)=\chi(G)-\cli(G)$;
\item
$K_n$ is the complete graph on $n$ vertices;
\item
$C_n$ is the simple cycle on $n$ vertices;
\item
$N_G(v)$ is the set of neighbours of a vertex $v$ in $G$.
\end{itemize}

Let $G_1$ and $G_2$ be two graphs. We denote by $G_1+G_2$ the graph $G$
for which $V(G)=V(G_1)\cup V(G_2)$, $E(G)=E(G_1)\cup E(G_2)\cup E'$,
where $E'=\{[x,y], x\in V(G_1),y\in V(G_2)\}$.

We will use the following theorem by \textsc{Dirac}~\cite{2}:

\begin{thm}\label{th:1.1}
Let $G$ be a graph such that $f(G)\ge 1$. Then $|V(G)|\ge\chi(G)+2$ and
$|V(G)|=\chi(G)+2$ only when $G=K_{\chi(G)-3}+C_5$.
\end{thm}

If $f(G)\ge 2$, then we have~\cite{12} (see also~\cite{16})

\begin{thm}\label{th:1.2}
Let $f(G)\ge 2$. Then
\begin{enumerate}[\indent\rm(a)]
\item
$|V(G)|\ge\chi(G)+4$;
\item
$|V(G)|=\chi(G)+4$ only when $\chi(G)\ge 6$ and $G=K_{\chi(G)-6}+C_5+C_5$.
\end{enumerate}
\end{thm}

In the case $\chi(G)=4$ and $\chi(G)=5$ we have the following more good
inequalities:
\begin{gather}\label{1.1}
\text{if $f(G)\ge 2$ and $\chi(G)=4$ then $|V(G)|\ge 11$, \cite{1};}\\
\text{if $f(G)\ge 2$ and $\chi(G)=5$ then $|V(G)|\ge 11$, \cite{13} (see also~\cite{14}).}
\label{1.2}
\end{gather}

For the case $f(G)\ge 3$ it was known that~\cite{12} (see also~\cite{17,18})

\begin{thm}\label{th:1.3}
Let $G$ be a graph such that $f(G)\ge 3$. Then $|V(G)|\ge\chi(G)+6$.
\end{thm}

In this paper we consider the case $|V(G)|=\chi(G)+6$. We prove the following
main theorem.

\begin{thm}\label{th:1.4}
Let $G$ be a graph such that $f(G)\ge 3$ and $|V(G)|=\chi(G)+6$.
Then $\chi(G)\ge 7$ and $G=K_{\chi(G)-7}+Q$ or $\chi(G)\ge 9$ and
$G=K_{\chi(G)-9}+C_5+C_5+C_5$, where $Q$ is the graph, whose
complementary graph $\overline Q$ is given in Fig.~1.
\end{thm}

\begin{figure}[t]
\centering
\includegraphics[width=50mm]{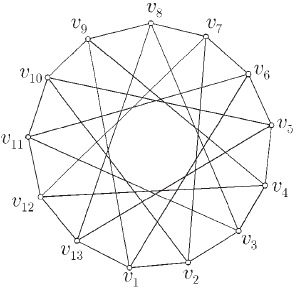}
\caption{}
\end{figure}

Obviously, if $f(G)\ge 3$ then $\chi(G)\ge 5$. Therefore we will consider only
the cases $\chi(G)\ge 5$. If $\chi(G)=5$ or $\chi(G)=6$ then by
Theorem~\ref{th:1.3} and Theorem~\ref{th:1.4} we see that $|V(G)|\ge\chi(G)+7$.
In these two cases we can state the following more strong results:
\begin{gather}\label{1.3}
\text{if $f(G)\ge 3$ and $\chi(G)=5$ then $|V(G)|\ge 22$, \cite{6};}\\
\text{if $f(G)\ge 3$ and $\chi(G)=6$ then $|V(G)|\ge 16$, \cite{9}.}
\label{1.4}
\end{gather}

The inequalities~\eqref{1.3} and~\eqref{1.4} are exact. \textsc{Lathrop} and
\textsc{Radziszowski}~\cite{9} proved that there are only two 16-vertex graphs
for which~\eqref{1.4} holds.

At the end of this paper we obtain by Theorem~\ref{th:1.4} new results about
some edge-Folkman numbers (Theorem~\ref{th:4.2}).

\section{Auxiliary results.}
A graph $G$ is defined to be vertex-critical chromatic if $\chi(G-v)<\chi(G)$
for all $v\in V(G)$. We shall use the following results of \textsc{Gallai}~\cite{4}
(see also~\cite{5}).

\begin{thm}\label{th:2.1}
Let $G$ be a vertex-critical chromatic graph and $\chi(G)\ge 2$.
If $|V(G)|<2\chi(G)-1$ then $G=G_1+G_2$, where $V(G_i)\ne\emptyset$, $i=1,2$.
\end{thm}

\begin{thm}\label{th:2.2}
Let $G$ be a vertex-critical $k$-chromatic graph, $|V(G)|$ and $k\ge 3$.
Then there exist $\ge\left\lceil\dfrac32\left(\dfrac53k-n\right)\right\rceil$
vertices with the property that each of them is adjacent to all the other
$n-1$ vertices.
\end{thm}

\begin{rem}\label{rem:2.1}
The formulations of Theorem~\ref{th:2.1} and Theorem~\ref{th:2.2} given above
are obviously equivalent to the original ones in~\cite{4} (see Remark~1 and
Remark~2 in~\cite{16}).
\end{rem}

\begin{pro}\label{pro:2.1}
Let $G$ be a graph such that $f(G)\ge 3$ and $|V(G)|=\chi(G)+6$. Then $G$ is
a vertex-critical chromatic graph.
\end{pro}

\begin{proof}
Assume the opposite. Then $\chi(G-v)=\chi(G)$ for some $v\in V(G)$.
Let $G'=G-v$. Since $\cli(G')\le\cli(G)$ we have $f(G')\ge f(G)\ge 3$.
By Theorem~\ref{th:1.3}
\[
|V(G')|\ge\chi(G')+6=\chi(G)+6=|V(G)|,
\]
which is a contradiction
\end{proof}

The following result by \textsc{Kerry}~\cite{7} will be used later.

\begin{thm}\label{th:2.3}
Let $G$ be a $13$-vertex graph such that $\alpha(G)\le 2$ and $\cli(G)\le 4$.
Then $G$ is isomorphic to the graph $Q$, whose complementary graph $\overline Q$
is given in Fig.~1.
\end{thm}

\begin{definition}\label{d:2.1}
The graph $G$ is called a Sperner graph if $N_G(u)\subseteq N_G(v)$ for
some $u,v\in V(G)$.
\end{definition}

Obviously if $N_G(u)\subseteq N_G(v)$ then $\chi(G-u)=\chi(G)$. Thus we have

\begin{pro}\label{pro:2.2}
Every vertex-critical chromatic graph is not a Sperner graph.
\end{pro}

The following lemmas are used in the proof of Theorem~\ref{th:1.4}.

\begin{lemma}\label{lem:2.1}
Let $G$ be a graph and $f(G)\ge 2$. Then
\begin{enumerate}[\indent\rm(a)]
\item
$|V(G)|\ge 10$;
\item
$|V(G)|=10$ only when $G=C_5+C_5$.
\end{enumerate}
\end{lemma}

\begin{proof}
The inequality~(a) follows from~\eqref{1.1}, \eqref{1.2} and
Theorem~\ref{th:1.2}(a). Let $|V(G)|=10$. Then by~\eqref{1.1},
\eqref{1.2} and Theorem~\ref{th:1.2}(a) we see that $\chi(G)=6$.
From Theorem~\ref{th:1.2}(b) we obtain $G=C_5+C_5$.
\end{proof}

\begin{lemma}\label{lem:2.2}
Let $G$ be a graph such that $f(G)\ge 3$ and $G$ is not a Sperner graph. Then
\[
|V(G)|\ge 11+\alpha(G).
\]
\end{lemma}

\begin{proof}
Assume the opposite, i.e.
\beq{2.1}
|V(G)|\le 10+\alpha(G).
\eeq
Let $A\subseteq V(G)$ be an independent set of vertices of $G$ such that
$|A|=\alpha(G)$. Consider the subgraph $G'=G-A$. From~\eqref{2.1} we see that
$|V(G')|\le 10$. Since $A$ is independent from $f(G)\ge 3$ it follows
$f(G')\ge 2$. According to Lemma~\ref{lem:2.1}(b), $G'=C_5^{(1)}+C_5^{(1)}$,
where $C_5^{(i)}$, $i=1,2$, are 5-cycles. Hence $\cli(G')=4$ and
$\cli(G)\le 5$. Thus if $a\in A$, then $N_G(a)\cap V(C_5^{(1)})$ or
$N_G(a)\cap V(C_5^{(2)})$ is an independent set.
Let $N_G(a)\cap V(C_5^{(1)})$ be independent set and
$C_5^{(1)}=v_1v_2v_3v_4v_5v_1$. Then we may assume that
$N_G(a)\cap V(C_5^{(1)})\subseteq\{v_1,v_3\}$. We obtain that
$N_G(a)\subseteq N_G(v_2)$ which contradicts the assumption of
Lemma~\ref{lem:2.2}.
\end{proof}

\begin{lemma}\label{lem:2.3}
Let $G$ be a graph such that $f(G)\ge 3$ and $|V(G)|=\chi(G)+6$.
Then $\chi(G)\ge 7$ and:
\begin{enumerate}[\indent\rm(a)]
\item
$G=Q$ if $\chi(G)=7$;
\item
$G=K_1+Q$ if $\chi(G)=8$;
\item
$G=K_2+Q$ or $G=C_5+C_5+C_5$ if $\chi(G)=9$.
\end{enumerate}
\end{lemma}

\begin{proof}
Since $\chi(G)\ne\cli(G)$ we have $\cli(G)\ge 2$. Thus, from $f(G)\ge 3$ it
follows $\chi(G)\ge 5$. By~\eqref{1.3} and~\eqref{1.4} we see that
$\chi(G)\ne 5$ and $\chi(G)\ne 6$. Hence $\chi(G)\ge 7$.

\textsc{Case~1.}
$\chi(G)=7$. In this case $|V(G)|=13$. From $\chi(G)=7$ and $f(G)\ge 3$
we see that $\cli(G)\le 4$. According to Proposition~\ref{pro:2.1} and
Proposition~\ref{pro:2.2}, $G$ is not a Sperner graph.
It follows from Lemma~\ref{lem:2.2} that $\alpha(G)\le 2$.
Thus, by Theorem~\ref{th:2.3}, $G=Q$.

\textsc{Case~2.}
$\chi(G)=8$. In this situation we have $|V(G)|=14$. By Proposition~\ref{pro:2.1},
$G$ is a vertex-critical chromatic graph. Since $|V(G)|<2\chi(G)-1$, from
Theorem~\ref{th:2.1} we obtain that $G=G_1+G_2$. Clearly,
\begin{align}
|V(G)|&=|V(G_1)|+|V(G_2)|;\label{2.2}\\
\chi(G)&=\chi(G_1)+\chi(G_2);\label{2.3}\\
f(G)&=f(G_1)+f(G_2);\label{2.4}\\[-7.25ex]\notag
\end{align}
\beq{2.5}
\text{$G_1$ and $G_2$ are vertex-critical chromatic graphs.}
\eeq

\textsc{Subcase 2.a.}
$G=K_1+G'$. Since $\chi(G')=7$ and $f(G')=f(G)\ge 3$, by the Case~1 we obtain
$G'=Q$ and $G=K_1+Q$.

\textsc{Subcase 2.b.}
$G_1$ and $G_2$ are not complete graphs. In this subcase, by~\eqref{2.5},
we have $\chi(G_i)\ge 3$ and $\chi(G_i)\ne\cli(G_i)$, $i=1,2$. Thus $f(G_i)\ge 1$,
$i=1,2$. According to Theorem~\ref{th:1.1}, $|V(G_i)|\ge 5$, $i=1,2$. From these
inequalities and~\eqref{2.2} it follows
\beq{2.6}
|V(G_i)|\le 9,
\quad
i=1,2.
\eeq
Let $f(G_1)\le f(G_2)$. Then, by~\eqref{2.4}, $f(G_2)\ge 2$.
From Lemma~\ref{2.1} we obtain $|V(G_2)|\ge 10$. This contradicts~\eqref{2.6}.

\textsc{Case 3.}
$\chi(G)=9$. In this case $|V(G)|=15$. By Proposition~\ref{pro:2.1}, $G$ is
a vertex-critical chromatic graph. Since $|V(G)|<2\chi(G)-1$,
from Theorem~\ref{th:2.1} it follows that $G=G_1+G_2$.

\textsc{Subcase 3.a.}
$G=K_1+G'$. Since $|V(G')|=14$, $\chi(G')=8$ and $f(G')=f(G)\ge 3$,
by Case~2 we have $G'=K_1+Q$. Hence $G=K_2+Q$.

\textsc{Subcase 3.b.}
$G_1$ and $G_2$ are not complete graphs. By~\eqref{2.5} it follows
$|V(G_i)|\ge 5$, $i=1,2$. From these inequalities and~\eqref{2.2}
we obtain
\beq{2.7}
|V(G_i)|\le 10,
\quad
i=1,2.
\eeq
Let $f(G_1)\le f(G_2)$. Then according to~\eqref{2.4} we have $f(G_2)\ge 2$.
From~\eqref{2.7} and Theorem~\ref{lem:2.1} we obtain $G_2=C_5+C_5$. Since
$|V(G_2)|=10$ and $\chi(G_2)=6$ we see from~\eqref{2.2} and~\eqref{2.3} that
$|V(G_1)|=5$ and $\chi(G_1)=3$. Thus, by~\eqref{2.5}, we conclude that
$G_1=C_5$. Hence $G=C_5+C_5+C_5$.
\end{proof}

\section{Proof of Theorem~\ref{th:1.4}.}
By Lemma~\ref{lem:2.3} we have that $\chi(G)\ge 7$. If $\chi(G)=7$ or
$\chi(G)=8$ Theorem~\ref{th:1.4} follows from Lemma~\ref{lem:2.3}.
Let $\chi(G)\ge 9$. We prove Theorem~\ref{th:1.4} by induction on $\chi(G)$.
The inductive base $\chi(G)=9$ follows from Lemma~\ref{lem:2.3}(c).
Let $\chi(G)\ge 10$. Then $\frac53\chi(G)-|V(G)|>0$.
By Proposition~\ref{pro:2.1} $G$ is vertex-critical chromatic graph.
Thus, according to Theorem~\ref{th:2.2}, we have $G=K_1+G'$.
As $\chi(G')=\chi(G)-1$, $f(G')=f(G)\ge 3$ and $|V(G')|=\chi(G')+6$,
we can now use the inductive assumption and obtain
\[
G'=K_{\chi(G')-7}+Q
\quad\text{or}\quad
G'=K_{\chi(G')-9}+C_5+C_5+C_5.
\]
Hence $G=K_{\chi(G)-7}+Q$ or $G=K_{\chi(G)-9}+C_5+C_5+C_5$.

\section{Edge Folkman numbers $F_e(a_1,\dots,a_r;R(a_1,\dots,a_r)-2)$.}
Let $a_1$, \dots, $a_r$ be integers, $a_i\ge 2$, $i=1,\dots,r$. The symbol
$G\toe(a_1\dots,a_r)$ means that in every $r$-coloring
\[
E(G)=E_1\cup\dots\cup E_r,
\quad
E_i\cap E_j=\emptyset,
\quad
i\ne j,
\]
of the edge set $E(G)$there exists a monochromatic $a_i$-clique $Q$ of colour
$i$ for some $i\in\{1,\dots,r\}$, that is $E(Q)\subseteq E_i$. The Ramsey
number $R(a_1,\dots,a_r)$ is defined as $\min\{n:K_n\toe(a_1,\dots,a_r)\}$.
Define
\begin{align*}
H_e(a_1,\dots,a_r;q)&=\{G:G\toe(a_1\dots,a_r)\text{ and }\cli(G)<q\};\\
F_e(a_1,\dots,a_r;q)&=\min\{|V(G)|:G\in H_e(a_1,\dots,a_r;q)\}.
\end{align*}
It is well known that
\beq{4.1}
F_e(a_1,\dots,a_r;q)\text{ exists }\iff q>\max\{a_1,\dots,a_r\}.
\eeq
In the case $r=2$ this was proved in~\cite{3} and the general case in~\cite{19}.
The numbers $F_e(a_1,\dots,a_r;q)$ are called edge Folkman numbers.
An exposition of the known edge Folkman numbers is given in~\cite{8}.
In this section we consider the numbers $F_e(a_1,\dots,a_r;R(a_1\dots,a_r)-2)$,
where $a_3\ge 3$, $i=1,\dots,r$. We know only one Folkman number of this kind,
namely $F_e(3,3,3;15)=23$ (see~\cite{11}).

In~\cite{12} we prove the following statement.

\begin{thm}\label{th:4.1}
Let $a_1,\dots,a_r$ be integers and $a_i\ge 3$, $i=1,\dots,r$, $r\ge 2$. Then
\beq{4.2}
F_e(a_1,\dots,a_r;R(a_1\dots,a_r)-2)\ge R(a_1\dots,a_r)+6.
\eeq
\end{thm}

\begin{rem}\label{rem:4.1}
It follows from $a_i\ge 3$ and $r\ge 2$ that $R(a_1,\dots,a_r)>2+
\max\{a_1,\dots,a_r\}$. Thus, by~\eqref{4.1}, the numbers
$F_e(a_1,\dots,a_r;R(a_1,\dots,a_r)-2)$ exist.
\end{rem}

The aim of this section is to prove the following result.

\begin{thm}\label{th:4.2}
Let $a_1,\dots,a_r$ be integers and $a_i\ge 3$, $i=1,\dots,r$, $r\ge 2$. Then
\[
F_e(a_1,\dots,a_r;R(a_1,\dots,a_r)-2)=R(a_1,\dots,a_r)+6
\]
if and only if $K_{R-7}+Q\toe(a_1,\dots,a_r)$ or $K_{R-9}+C_5+C_5+C_5\toe
(a_1,\dots,a_r)$, where $R=R(a_1,\dots,a_r)$.
\end{thm}

We shall use the following result obtained by \textsc{Lin}~\cite{10}:
\beq{4.3}
G\toe(a_1,\dots,a_r)
\Rightarrow
\chi(G)\ge R(a_1,\dots,a_r).
\eeq

\begin{proof}[Proof of Theorem~\ref{th:4.2}]
I. Let $F_e(a_1,\dots,a_r;R-2)=R+6$. Let $G\in H_e(a_1,\dots,a_r;R-2)$ and
\beq{4.4}
|V(G)|=R+6.
\eeq
Since $\cli(G)\le R-3$, from~\eqref{4.3} it follows $f(G)\ge 3$.
By Theorem~\ref{th:1.3}, we have
\beq{4.5}
|V(G)|\ge\chi(G)+6.
\eeq
From~\eqref{4.3}, \eqref{4.4} and~\eqref{4.5} we see that $\chi(G)=R$ and
$|V(G)|=\chi(G)+6$. Thus, according to Theorem~\ref{th:1.4}, $G=K_{\chi(G)-7}+Q=
K_{R-7}+Q$ or $G=K_{\chi(G)-9}+C_5+C_5+C_5=K_{R-9}+C_5+C_5+C_5$. This implies
$K_{R-7}+Q\toe(a_1,\dots,a_r)$ or $K_{R-9}+C_5+C_5+C_5\toe(a_1,\dots,a_r)$
because $G\in H_e(a_1,\dots,a_r;R-2)$.

II. Let $K_{R-7}+Q\toe(a_1,\dots,a_r)$. Then $K_{R-7}+Q\in H_e(a_1,\dots,a_r;R-2)$
because $\cli(K_{R-7}+Q)=R-3$. Hence
\[
F_e(a_1,\dots,a_r;R-2)\le|V(K_{R-7}+Q)|=R+6.
\]
This inequality and~\eqref{4.2} imply that $F_e(a_1,\dots,a_r;R-2)=R+6$.

In the same way we see that from $K_{R-9}+C_5+C_5+C_5\toe(a_1,\dots,a_r)$
it follows that $F_e(a_1,\dots,a_r;R-2)=R+6$.
\end{proof}

\begin{rem}\label{rem:4.2}
We obtain the equality $F_e(3,3,3;15)=23$ proving that
$K_8+C_5+C_5+C_5\toe(3,3,3)$. We do not know whether
$K_{10}+Q\toe(3,3,3)$.
\end{rem}

\begin{rem}\label{rem:4.3}
By Theorem~\ref{th:4.1} we have $F_e(3,5;12)\ge 20$ and $F_e(4,4;16)\ge 24$.
The exact values of these numbers are not known. Therefore, having in mind
Theorem~\ref{th:4.2}, it will be interesting to know whether the following
statements are true:
\begin{align*}
&K_7+Q\toe(3,5),&&K_5+C_5+C_5+C_5\toe(3,5);\\
&K_{11}+Q\toe(4,4),&&K_9+C_5+C_5+C_5\toe(4,4).\\
\end{align*}
\end{rem}

\begin{rem}\label{rem:4.4}
By Theorem~\ref{th:4.1}, $F_e(3,4;7)\ge 15$. It was proved in~\cite{8} that
$F_e(3,4;8)=16$. Thus $F_e(3,4;7)\ge 17$.
\end{rem}

\oneaddr{%
Faculty of Mathematics and Informatics\\
St Kliment Ohridski University of Sofia\\
5, James Bourchier Blvd\\
1164 Sofia, Bulgaria\\
\email{nenov@fmi.uni-sofia.bg}}

\end{document}